\documentclass[titlepage,11pt]{article}
\usepackage{latexsym}
\usepackage{amsfonts}
\usepackage{amsmath}
\usepackage{amssymb}
\usepackage{mathtools}
\usepackage{tikz}
\usepackage{comment}
\newcommand\blackslug{\hbox{\hskip 1pt \vrule width 4pt height 8pt depth 1.5pt
		\hskip 1pt}}
\newcommand\bbox{\hfill \quad \blackslug \bigbreak}

\def\LL{,\ldots,}

\newcommand{\mac}{\mathcal}

\title{Line-width and path-width
}
\author{
	Tung Nguyen\\
	University of Oxford, \\
	Oxford, UK\\
	\and
	Alex Scott\\
	University of Oxford, \\
	Oxford, UK
	\and
	Paul Seymour\\
	Princeton University,\\ Princeton, NJ 08544, USA}

\date{
}

\newtheorem{thm}{}[section]

\newcommand{\Proof}{\noindent{\bf Proof.}\ \ }

\begin{document}
	\maketitle
	\begin{abstract}
		For finite graphs, path-width is an interesting and useful concept, but if we extend it to infinite graphs in the most obvious way (by making the indexing path infinite), it does not work nicely.
		The simplest extension that works nicely is to allow the indexing set to be any totally-ordered set, and then the 
		corresponding decomposition is called a ``line-decomposition'', and the maximum possible width is called ``line-width''.  
		
		In particular, the indexing set need not be a well-order; but the corresponding decomposition would be easier to use if it 
		was. We show that if a graph has line-width at most $k$, it admits a well-ordered line-decomposition with width at most $2k$, 
		and this is best possible.

	\end{abstract}
	
	\section{Introduction}
	
	A {\em tree-decomposition} of a finite graph $G$ is a pair $(T,(W_t:t\in V(T)))$, where $T$ is a finite tree, and $W_t$
	is a subset of $V(G)$ for each $t\in V(T)$ (called a {\em bag}), such that:
	\begin{itemize}
		\item $G=\bigcup_{t\in V(T)} G[W_t]$, and 
		\item for all $t_1,t_2,t_3\in V(T)$, if $t_2$ lies on the path of $T$ between $t_1,t_3$, then
		$B_{t_1}\cap B_{t_3}\subseteq B_{t_2}$.
	\end{itemize}
	The {\em width} of a tree-decomposition $(T,(W_t:t\in V(T)))$ is the maximum of the numbers $|W_t|-1$ for $t\in V(T)$,
	and the {\em tree-width} of $G$ is the minimum width of a tree-decomposition of $G$.
	
	Similarly, a 
	{\em path-decomposition} of a finite graph $G$ is the same except that we require $T$ to be a path; and
	the {\em path-width} of a finite graph is the minimum width
	of a path-decomposition. 
	
	We extend tree-width and path-width to infinite graphs by 
	allowing the indexing tree or path to be infinite
	(thus, for a path-decomposition, the indexing path might be finite, or one-way infinite, or two-way infinite).
	This seems the most natural extension, and for tree-width it works nicely. In particular, there is a ``compactness'' theorem of
	Thomas~\cite{thomas} (and see also~\cite{subtrees} for a proof): 
	\begin{thm}\label{thomasthm}
		A graph $G$ has tree-width at most $k$ if and only if every finite subgraph of $G$ has tree-width at most $k$.
	\end{thm}
	There are several ways to define infinite analogues of finite trees: for instance, topological trees, or set-theoretic trees.
	We could choose our favourite one of these definitions (let us call them ``supertrees''), 
	and use a supertree as the indexing set,  giving ``supertree-decompositions'' and ``supertree-width''.
	But this gives nothing new: any graph with supertree-width at most $k$ also has tree-width at most $k$. To see this,
	observe that, whatever definition of supertree we have taken, if $G$ has supertree-width at most $k$, then
	all its finite subgraphs have tree-width at most $k$, and therefore $G$ has tree-width at most $k$ by Thomas' theorem \ref{thomasthm}.

	For path-width, the analogue of Thomas' theorem is far from true. For instance, all graphs with finite path-width have only 
	countably many vertices, so a graph with uncountably many vertices and no edges is already a counterexample; all its finite subgraphs have path-width zero, and yet $G$ does not have finite path-width. Even for countable graphs, the analogue is still not true: for instance,
	any graph with infinitely many infinite connected subgraphs that are pairwise vertex-disjoint does not have finite path-width,
	and yet its finite subgraphs could all have very small path-width. 
	
	Thus, for an infinite graph $G$,  all finite subgraphs having path-width at most
	$k$ does not imply that $G$ has path-width at most $k$. But what structure of $G$ does it imply?
	This was answered in~\cite{subtrees} as follows.

	A {\em line} is a nonempty set $L$ together with a linear order $\le_L $.  A {\em line-decomposition} of $G$ is a pair
	$(L,(W_t:t\in L))$ where $L$ is a line, and each $W_t$ is a subset of $V(G)$, satisfying:
	\begin{itemize}
		\item $G=\bigcup_{t\in L} G[W_t]$, and 
		\item $W_t\cap W_{t''}\subseteq W_{t'}$ for all $t,t',t''\in L$ with $t\le_L t'\le_L t''$.
	\end{itemize}
	The second condition here is called the ``betweenness axiom''.
	We define the {\em width} of such a decomposition to be the maximum of $|B_t|-1$ over all $t\in L$, if this maximum exists, and $\infty$ otherwise; and the {\em line-width} of $G$ is the minimum width of a line-decomposition.
	(This definition was used in~\cite{AS2, AS3}, to increase the generality of some results about graphs with bounded path-width.
	It is also used in~\cite{diestel, thomas2}, where what we call ``line-decompositions'' are called ``linear decompositions''.)
	
	With Chudnovsky, we proved in~\cite{subtrees} that:
	\begin{thm}\label{linecompact}
		A graph $G$ has line-width at most $k$ if and only if every finite subgraph of $G$ has path-width at most $k$.
	\end{thm}
	But do we need such general linear orders? Could we restrict the order-type somehow, perhaps in terms of the 
	maximum bag-size? Could we make the order a well-order? This paper is an investigation into these and related questions.
	In particular, we will show that:
	\begin{thm}\label{mainthm}
		If a graph $G$ has finite line-width $k$, then it admits a line-decomposition $(L,(W_t:t\in L))$ of width at most $2k$
		such that $L$ is a well-order.
	\end{thm}
	We will also show that $2k$ is best possible. 
	
	There is a very similar result of Erd\H{o}s and Hajnal~\cite{hajnal}: for a graph $G$, if every finite subgraph admits a linear order of its vertex set 
	such that each vertex is adjacent to at most $k$ earlier vertices, then there is a well-order of $V(G)$ such that each vertex is adjacent to at most 
	$2k$ earlier vertices. Again $2k$ is best possible.
	
	At the end of the paper, we look at a different question. 
	We will show that, after a little tidying,  every line-decomposition of width $\le k$ of a graph $G$ can be built by substitution and concatenation,
	starting from so-called ``prime'' line-decompositions, which are much simpler.
	
	\section{Tidying a line-decomposition}
	
	If $L$ is a line, an {\em interval} of $L$ is a nonempty subset $I\subseteq L$ such that if $r,s,t\in L$ with $r\le_L s\le_L t$, and $r,t\in I$,
	then $s\in I$. 
	An interval $I$ is {\em initial} if $I\ne L$ and
	for all $s,t\in L$ with $s\le_L t$, if $t\in I$ then $s\in I$. (Note that if $I$ is an initial interval then $I, L\setminus I$ are 
	both nonempty, and so $L\setminus I$ is an initial interval of the line obtained from $L$ by reversing its order.) 
	
	If $(L,(W_t:t\in L))$ is a line-decomposition of $G$, and $I$ is an interval, we define
	$$W(I) = \bigcup_{t\in I}W_t.$$
	If $I$ is an initial interval, then $W(I)\cap  W(L\setminus I)$
	is called the {\em $I$-split}; every path between $W(I),  W(L\setminus I)$ has a vertex in this set. 
	We say $X$ is a {\em split} if $X$ is the $I$-split for some initial interval $I$.
	
	Let us see first some easy things we can do to make a line-decomposition nicer without increasing its width. 
	Let $(L,(W_t:t\in L))$ be a line-decomposition of $G$, of finite width. If there are distinct $s,s'\in L$ with $W_s=W_{s'}$,
	we could remove one of $s,s'$ from $L$ and still have a line-decomposition. In general, we can partition 
	$L$ into equivalence
	classes under the relation $s\equiv s'$ if $W_s=W_{s'}$, and choose one element of each equivalence class. Let $L'$ be the set of
	chosen elements, ordered by the order inherited from $L$. Then $(L',(W_t:t\in L'))$ is still a  line-decomposition of $G$, with 
	the same width, and now there are no repeated bags.
	So we could assume that all the sets $W_t\;(t\in L)$ are distinct. 
	
	Now let $L'$ be the set of all $s\in L$ such that there is no $s'\in L$ with $s'\ne s$ and with $W_{s}\subseteq W_{s'}$;
	then $(L',(W_t:t\in L'))$ is also a line-decomposition of $G$, with the same width. So we could assume that no $W_s$ is a subset of another.
	
	Let us add edges to $G$, forming $G'$, in such a way that every edge we add has both ends in some $W_t$, and every $W_t$ becomes a clique. The same
	line-decomposition is also a line-decomposition of $G'$, and any line-decomposition
	of the new graph is also a line-decomposition of $G$; so we might as well work with the new graph in place of $G$. 
	
	Thus, 
	so far we are assuming:
	\begin{itemize}
		\item $G$ is non-null;
		\item $(L,(W_t:t\in L))$ is a line-decomposition of finite width;
		\item for all distinct $s,t\in L$, $W_s\not\subseteq W_t$; and
		\item $W_t$ is a clique of $G$ for each $t\in L$.
	\end{itemize}
	Let us say a line-decomposition  $(L,(W_t:t\in L))$ is {\em tidy} if it satisfies these conditions. If $(L,(W_t:t\in L))$ is tidy, then 
	each $W_t$ is a clique, and every clique of $G$ is a subset of some $W_t$ (this is an easy exercise); and consequently, each $W_t$ is 
	a maximal clique of $G$, and every maximal clique of $G$ is one of the sets $W_t$. Thus,  the sets $W_t$ are precisely the maximal cliques of $G$, and in particular they are determined by $G$; 
	the only freedom in choosing a tidy line-decomposition is choosing the order  $\le_LL$, that is, choosing a linear order of the set of maximal cliques 
	that satisfies the betweenness axiom. (See~\cite{subtrees, tucker} for when such an order exists.) Thus, the cardinality of $L$ 
	is the cardinality of the set of maximal cliques (that is, $|G|$ if $G$ is infinite).
	There can be many different such orders; for instance, if $G$ is a ``star'' graph, a tree in which all 
	edges have a common end, 
	then the maximal cliques of $G$ are just the edges, and can be arranged in any linear order to make a line-decomposition.
	
	This tidying process is identical to  a process of Halin~\cite{halin}, and the corresponding line-decomposition is called by Halin a ``consecutive generalized prime graph decomposition''. 
	
	We observe: 
	\begin{thm}\label{splitlemma}
		Let $(L,(W_t:t\in L))$ be a line-decomposition of $G$, and let $I\subseteq L$ be an initial interval such that the split $S=W(I)\cap W(L\setminus I)$ is 
		finite. Then there exist $i\in I$ and $j\in L\setminus I$ such that $S=W_i\cap W_j$.
	\end{thm}
	\Proof
	If $S=\emptyset$ then any $i\in I$ and any $j\in L\setminus I$ would suffice; so we may assume that $S\ne\emptyset$.
	Then for each $v\in S$, choose $i_v\in I$ and $j_v\in L\setminus I$ with $v\in W_{i_v}\cap W_{j_v}$. Let $i = \max(i_v:v\in S)$
	and $j = \min(j_v:v\in S)$, both of which exist since $S$ is finite.  Then $i\le_L j$, and the betweenness axiom implies that 
	$S\subseteq W_i\cap W_j$. The reverse inclusion is immediate. This proves \ref{splitlemma}.~\bbox
	
	Consequently, each finite split is contained in two bags. If the line-decomposition is tidy, each split is therefore a proper subset of a bag (since no two bags are equal), and therefore has size less than the maximum bag size.
	
	To aid the reader's intuition, let us mnetion an observation of a referee,  that if $G$ is connected and admits a tidy line-decomposition of finite width, then
	the splits of the decomposition are precisely the minimal finite cutsets of $G$. We omit the proof since we do not need the result.

	\section{Well-ordered line-decompositions}
	
	A line $L$ is a {\em well-order} if there is no infinite sequence $t_1,t_2,\ldots  $ of distinct elements of $L$
	such that $t_{i+1}\le_L t_i$ for each $i\ge 1$; and a line-decomposition $(L,(W_t:t\in L))$ is a {\em wo-decomposition}
	if $L$ is a well-order. The minimum width of a wo-decomposition of $G$ is called the {\em wo-width} of $G$.
	
	Our objective is to show that the wo-width of $G$ is at most twice the line-width of $G$. Let us see first that this is best possible. 
	\begin{thm}\label{bestposs}
		For each integer $k\ge 1$, there is a graph with line-width at most $k$ and with wo-width at least $2k$.
	\end{thm}
	\Proof Let $G$ be the graph with vertex set $\mathbb Z$, where for all distinct integers $i,j$,  we say $i,j$ are 
	adjacent if $|j-i|\le k$. Let $B_i=\{i,i+1\LL i+k\}$ for each $i\in \mathbb Z$; then $(B_i:i\in \mathbb{Z})$ is a line-decomposition
	of $G$ of width $k$. Now suppose that $(L, (W_t:t\in L))$ is a wo-decomposition of $G$ of width less than $2k$, and so $|W_t|\le 2k$ for each $t\in L$. 
	Every clique of $G$
	is a subset of some $W_t$. Since $L$ is a well-order, there exists $t\in L$ minimum such that $W_t$ includes $B_i$ for some integer $i$.
	Choose $t', t''\in L$ such that $B_{i+k+1}\subseteq W_{t'}$ and $B_{i-k-1}\subseteq W_{t''}$. Thus $t\le_L t',t''$ from the choice of $t$. 
	Since $B_i, B_{i+k+1}, B_{i-k-1}$ are pairwise disjoint and have size $k+1$, no $W_t$ includes two of them, since each $W_t$ has size at most $2k$. In particular,
	$t, t', t''$ are all different. There is a matching $M$ of size $k$ between 
	$B_i, B_{i+k+1}$, joining $j$ and $j+k$ for $i+1\le j\le i+k$. Since $|W_{t''}|\le 2k$, and $B_{i-k-1}\subseteq W_{t''}$, there is 
	an edge of $M$ with neither end in $W_{t''}$; and so it is not the case that $t\le_L t'' \le_L t'$. But similarly, there
	is a matching of size $k$ between $B_i, B_{i-k-1}$, and so it is not the case that $t\le_L t' \le_L t''$, a contradiction.
	This proves \ref{bestposs}.~\bbox
	
	Let $L$ be a line. We say that $L$ is {\em integral} if $L$ has order-type that of a subset of the set of integers, that is, 
	there is a function $\phi:L\to\mathbb{Z}$ such that
	if $s,t\in L$ are distinct and $s\le_L t$ then $\phi(s)<\phi(t)$.
	\begin{thm}\label{trivialorder}
		Let $L$ be a line, and suppose that for all $r,t\in L$ with $r\le_L t$, there are only finitely many $s\in L$ with $r\le_L s\le_L t$.
		Then $L$ is integral.
	\end{thm}
	\Proof
	Choose $r\in L$ and define $\phi(r)=0$. For each $t\in L\setminus \{r\}$ with $r<_Lt$, let $\phi(t)$ be the number of $s\in L$ with
	$r<_Ls\le_L t$. If $t<_Lr$, let $-\phi(t)$ be  the number of $s\in L$ with
	$t\le_Ls<_Lr$. Then $\phi$ satisfies the theorem. This proves \ref{trivialorder}.~\bbox

	For a split $S$, there may be many different initial intervals $I$ such that $S$ is the $I$-split. For instance, when $G$
	is the star graph, and $v\in V(G)$ is an end of every edge of $G$, then $\{v\}$ is the $I$-split for every initial interval $I$, for every tidy line-decomposition.
	
	Let $(L,(W_t:t\in L))$ be a line-decomposition of $G$, and let $S_1,S_2$ be splits.  We say $S_1$ is {\em before} $S_2$ 
	if one the following holds:
	\begin{itemize}
		\item $S_1=S_2$; or
		
		\item $S_1\ne S_2$, and $I_1\subseteq I_2$ for all choices of initial intervals $I_1,I_2$ such that $S_1$, $S_2$ are the $I_1$-split and $I_2$-split respectively.
	\end{itemize} 
	
	\begin{thm}\label{ordered}
		Let $(L,(W_t:t\in L))$ be a line-decomposition of $G$, and let $S_1,S_2$ be distinct splits of the same finite size. Then either $S_1$ is before $S_2$, or $S_2$ is before $S_1$.
	\end{thm}
	\Proof For $i = 1,2$, let $\mac I_i$ be the set of all initial intervals $I$ such that $S_i$ is the $I$-split. 
	\\
	\\
	(1) {\em For all $I_1\in \mac I_1$ and $I_2\in \mac I_2$, if $I_1\subseteq I_2$ then $I_1\subseteq I_2'$ for all $I_2'\in \mac I_2$, and $I_1'\subseteq I_2$ for all $I_1'\in \mac I_1$.}
	\\
	\\
	Suppose that there exists $I_2'\in \mac I_2$ such that $I_1\not\subseteq I_2'$. Thus, $I_2'\subseteq I_1$.  Since $S_1\ne S_2$
	and $|S_1|=|S_2|$ is finite, there exists $v\in S_2\setminus S_1$. Since $v\in S_2$ and $S_2$ is the $I_2$-split, it follows that 
	$v\in W(L\setminus I_2)\subseteq W(L\setminus I_1)$; and similarly, since $S_2$ is the $I_2'$-split, $v\in W(I_2')\subseteq W(I_1)$. But then $v\in W(I_1)\cap W(L\setminus I_1)=S_1$, a contradiction. This proves the first statement of (1), and the second is 
	proved similarly. This proves (1). 
	
	\bigskip
	
	Choose $I_1\in \mac I_1$ and $I_2\in \mac I_2$; then one of $I_1,I_2$ is a subset of the other, and by exchanging $S_1,S_2$
	if necessary we may assume that $I_1\subseteq I_2$. Now let $I_1'\in \mac I_1$ and $I_2'\in \mac I_2$.
	By the first statement of (1), we have $I_1\subseteq I_2'$, and so by the second statement of (1), with $I_2$ replaced by $I_2'$,
	we have 
	$I_1'\subseteq I_2'$.  This proves that $S_1$ is before $S_2$, and so proves \ref{ordered}.~\bbox
	
	As we said earlier, one split can be the $I$-split for many initial intervals $I$; but we can enumerate the {\em distinct} 
	splits of minimum size by the following:
	\begin{thm}\label{countable}
		Let $(L,(W_t:t\in L))$ be a tidy line-decomposition of $G$, such that for some $m$,
		every split has size at least $m$. Then we can number the splits with size exactly $m$ as
		$\{S_i:i\in K\}$, where $K$ is an interval of $\mathbb{Z}$, in such a way that $S_i$ is before $S_j$ for $i,j\in K$ with $i< j$.
	\end{thm}
	\Proof
	Let $\mac S$ be the set of all splits with size $m$. By \ref{ordered}, the relation ``is before'' is a linear order on $\mac S$, and we need to show it is integral. For $S,S'\in \mac S$, we write $S\le S'$ if $S$ is before $S'$. 
	Suppose that there exist $S_1,S_2\in \mac S$ such that $S_1\le S_2$, and there are infinitely many $S\in \mac S$ 
	such that $S_1\le S\le S_2$.
	Choose a subset $X\subseteq V(G)$, maximal such that there are infinitely many $S\in \mac S$ with
	$S_1\le S\le S_2$ and with $X\subseteq S$. (This is possible since every such $X$ has size at most the width.)
	A path of $G$ {\em avoids} $X$ if it has no vertex in $X$.
	\\
	\\
	(1) {\em There is no path of $G\setminus X$ between $S_1, S_2$
		avoiding $X$.}
	\\
	\\
	Suppose that there is a path $P$ of $G$ between $S_1, S_2$
	avoiding $X$,
	and choose such a path $P$ minimal. Thus, none of its internal vertices belong to $S_1\cup S_2$.
	But for each $S\in \mac S$ with $S_1\le S\le S_2$,
	some vertex of $P$ belongs to $S$; and since $P$ has finite length, there is a vertex $v\in V(P)$ such that $v$
	belongs to $S$
	for infinitely many $S$ with $S_1\le S\le S_2$, contradicting the maximality of $X$. This proves (1).
	
	\bigskip
	
	For $i = 1,2$, choose an initial interval $I_i$ such that $S_i$ is the $I_i$-split,
	and let $I$ be the set of all $t\in 
	I_2$ such that either $t\in I_1$ or there is a path of $G$ between $S_1,W_t$ avoiding $X$.
	\\
	\\
	(2) {\em $I$ is an initial interval}.
	\\
	\\
	Let  $s,t\in L$ with $s\le_L t$ and $t\in I$; we need to show that $s\in I$. If $s\in I_1$ the claim is true, so we assume that
	$s\notin I_1$;
	and so $s,t\in I\setminus I_1$. Hence there is a path of $G$ between $S_1,W_t$ avoiding $X$; but this path intersects $W_s$,
	and so $s\in I$. This proves (2).
	
	\bigskip
	
	Let $S$ be the $I$-split.
	\\
	\\
	(3) {\em $S\subseteq X$.}
	\\
	\\
	Suppose
	that there exists $v\in S\setminus X$. Since $v\in W(L\setminus I)$, there exists $t\in L\setminus I$ with $v\in W_t$; and hence $t\notin I_1$ and
	there is no path of $G$ between $S_1,W_t$ avoiding $X$.  In particular, $v\notin S_1$.
	Since $v\in W(I)$, there exists $s\in I$ with $v\in W_s$; and therefore either $s\in I_1$ or
	there is a path of $G$ between $S_1,W_s$ avoiding $X$. Since $v\in W_s\cap W_t$ and $v\notin S_1$, it follows that
	$s\notin I_1$. Hence there is a path of $G$ between $S_1, W_s$ avoiding $X$, and so there is a path of $G$ between $S_1,v$
	avoiding $X$, since $v\in W_s\setminus X$ and $W_s$ is a clique. But there is no path of $G$ between $S_1, W_t$
	avoiding $X$, a contradiction. This proves (3).
	
	\bigskip
	
	Since $S$ is a split, it follows from the hypothesis that
	$|S|\ge m$. But $X$ is a subset of infinitely many distinct splits all of size $m$, and so $|X|<m$, contrary to (3). This proves \ref{countable}.~\bbox

	A referee observed that \ref{countable} (and step (3) in the proof of \ref{addlimit} below) can also be proved using a result of Halin~\cite[Corollary 1]{halin2} that for every integer $k\ge1$ and every graph $G$ with distinct vertices $u,v$, there are only finitely many minimal $u$-$v$ separators of size at most $k$ in $G$.
	However, we decided to present our current proofs in full to keep the paper self-contained.

	If $(L,(W_t:t\in L))$ is a line-decomposition of $G$, we say a vertex $v$ of $G$ is a {\em left-limit vertex} of the line-decomposition
	if for all $s,t\in L$ with $s\le_L t$, if $v\in W_t$ then $v\in W_s$. (There may be no such vertices.) Similarly, $v$ is a {\em right-limit vertex} if
	for all $s,t\in L$ with $s\le_L t$, if $v\in W_s$ then $v\in W_t$.
	\begin{thm}\label{doublesplit}
		Let $(L,(W_t:t\in L))$ be a tidy line-decomposition of $G$,  such that for some $m$, 
		every split has size at least $m$. Let $S$ be a split of size $m$, and let $\mac I$ be the set of all initial intervals $I$ such 
		that $S$ is the $I$-split. Then either $\bigcap_{I\in \mac I}I=\emptyset$ and $S$ is the set of all left-limit vertices of 
		$(L,(W_t:t\in L))$, or $\bigcap_{I\in \mac I}I$ is the minimal 
		member of $\mac I$. Similarly, either $\bigcup_{I\in \mac I}I=L$ and $S$ is the set of all right-limit vertices, or $\bigcup_{I\in \mac I}I$ is the maximal member of $\mac I$. 
	\end{thm}
	\Proof
	Let $I_1=\bigcap_{I\in \mac I}I$. Suppose first that $I_1=\emptyset$. We must show that $S$ is the set of all left-limit vertices. 
	We show first that every vertex in $S$ is a left-limit vertex. Let $v\in S$ and let $s,t\in L$ with $s\le_L t$ and $v\in W_t$. 
	Since $\bigcap_{I\in \mac I}I=\emptyset$, there exists $I\in \mac I$ with $s\notin I$ and hence $t\notin I$. Since $S$ is the 
	$I$-split, $v\in W(I)$, and so $v\in W_r$ for some $r\in I$; but $r\le_Ls\le_L t$ and $v\in W_r\cap W_t$, and therefore $v\in W_s$. 
	This proves that every vertex in $S$ is a left-limit vertex. For the converse inclusion, we must show that if $I_1=\emptyset$
	then every left-limit vertex belongs to $S$. Thus, let $v$ be a left-limit vertex, and choose $t\in L$ with $v\in W_t$.
	Since $\bigcap_{I\in \mac I}I=\emptyset$, there exists $I\in \mac I$ with $t\in L\setminus I$. Choose $s\in I$; then $s\le_L t$,
	and so $v\in W_s$ since $v$ is a left-limit vertex; and consequently $v\in W(I)\cap W(L\setminus I) = S$. This proves that if
	$\bigcap_{I\in \mac I}I=\emptyset$ then $S$ is the set of all left-limit vertices.
	
	Next we assume that $I_1\ne \emptyset$.
	Since $\mac I\ne \emptyset$, it follows that $I_1$ is an initial interval; let $S_1$ be the $I_1$-split. If $S_1\subseteq S$,
	then since $S_1$ is a split, and $S$ is a split of minimum size, it follows that $S_1=S$ and therefore $I_1$ 
	is the minimal          
	member of $\mac I$; so we assume for a contradiction that there exists $v\in S_1\setminus S$. Since $v\in S_1$, there exists
	$s\in I_1$ and $t\in L\setminus I_1$  such that $v\in W_s\cap W_t$. Since $t\notin I_1$, 
	there exists $I\in \mac I$ with $t\notin I$ and therefore $t\in L\setminus I$. But $s\in I_1\subseteq I$, and 
	$$v\in W_s\cap W_t\subseteq W(I)\cap W(L\setminus I) = S,$$
	a contradiction. This proves the first statement, and the second follows similarly. 
	
	This proves \ref{doublesplit}.~\bbox

	We say $(L,(W_t:t\in L))$ is a
	line-decomposition of $G$ {\em from $Z_1$ to $Z_2$} if $Z_1$ is a set of left-limit vertices and $Z_2$ is a set of right-limit
	vertices. (We emphasize that $Z_1$ need not be the entire set of left-limit vertices, just a subset, and similarly for $Z_2$.) We will need the following lemma:
	\begin{thm}\label{concat1}
		Let $(L,(W_t:t\in L))$ be a line-decomposition of $G$ from $Z_1$ to $Z_2$, of width at most $k$, and let $I$ be an initial interval.
		Let $S$ be the $I$-split, and let $G_1=G[W(I)]$ and $G_2 = G[W(L\setminus I)]$. Then:
		$(I, (W_t:t\in I))$ is a line-decomposition of $G_1$ from $Z_1$ to $S$,  and $(L\setminus I,(W_t:t\in L\setminus I))$ is a 
		line-decomposition of $G_2$ from $S$ to $Z_2$, both of width at most $k$.
	\end{thm}
	\Proof
	Certainly $(I, (W_t:t\in I))$ is a line-decomposition of $G_1$, and we need to check it is from 
	$Z_1$ to $S$. Thus, we need to check that every member of $Z_1$ is a left-limit vertex of $(I, (W_t:t\in I))$, and every vertex 
	in $S$ is a right-limit vertex of $(I, (W_t:t\in I))$. 
	
	Let $v\in Z_1$. Then $v\in W_t$ for some $t\in L$, and $v\in W_s$
	for all $s\in L$ with $s\le_L t$. Since $I\ne \emptyset$, there exists $s\in I$, and hence either $t\in I$ or $v\in W_s$, and in the second case we may replace $t$ by $s$. Thus, 
	we may assume that $t\in I$. Since $v\in W_r$
	for all $r\in L$ with $r\le_L t$, it follows that $v$ is a left-limit vertex of $(I, (W_t:t\in I))$. So 
	every member of $Z_1$ is a left-limit vertex of $(I, (W_t:t\in I))$. 
	
	Now let $u\in S$. Since 
	$S$ is the $I$-split, there exists $s\in I$ and $t\in L\setminus I$ such that $u\in W_s\cap W_t$; but then $u\in W_{s'}$
	for all $s'\in I$ with $s\le_L s'$, since $s\le_L s'\le_L t$. Consequently $u$ is a right-limit vertex of $(I, (W_t:t\in I))$,
	and so every vertex
	in $S$ is a right-limit vertex of $(I, (W_t:t\in I))$.
	
	This proves that $(I, (W_t:t\in I))$ is a line-decomposition of $G_1$ from $Z_1$ to $S$. Similarly, 
	$(L\setminus I,(W_t:t\in L\setminus I))$ is a 
	line-decomposition of $G_2$ from $S$ to $Z_2$, and clearly both these decompositions have width at most $k$. This proves 
	\ref{concat1}.~\bbox
	
	\begin{thm}\label{concat2}
		Let $Z_1,Z_2,S\subseteq V(G)$, and let $G_1,G_2$ be induced subgraphs of $G$ with $G_1\cup G_2 = G$ and $V(G_1)\cap V(G_2) = S$,
		and $Z_i\subseteq V(G_i)$ for $i = 1,2$. 
		If $G_1$ admits a wo-decomposition $(L_1,(W_t:t\in L_1))$ from $Z_1$ to $S$,
		and $G_2$ admits a wo-decomposition $(L_2,(W_t:t\in L_2))$ from $S$ to $Z_2$, both of width at most some $k'$, then concatenating them yields a wo-decomposition of $G$ from
		$Z_1$ to $Z_2$ of width at most $k'$ and of order type $L_1+L_2$.
	\end{thm}
	\Proof 
	We may assume that $L_1\cap L_2=\emptyset$. Let $L$ be the linear order on $L_1\cup L_2$ defined by
	$x\le_L y$ if and only if either $x\in L_1$ and $y\in L_2$, or for some $i\in \{1,2\}$, $x,y\in L_i$ and $x\le_{L_i}y$. 
	To see that $(L,(W_t:t\in L))$ is a line-decomposition of $G$, observe first that 
	$$\bigcup_{t\in L}G[W_t] = \bigcup_{t\in L_1}G[W_t]\cup \bigcup_{t\in L_2}G[W_t] = G_1\cup G_2 = G.$$
	Now suppose that $r\le_L s\le_L t$; we must show that $W_r\cap W_t\subseteq W_s$. This is true if $t\in L_1$ (because then $r,s\in L_1$), or if $r\in L_2$; so we assume that $r\in L_1$ and $t\in L_2$. Now $s$ may belong to either of $L_1,L_2$, but from the symmetry we may assume that $s\in L_1$. (This is a small cheat, since we are working with well-orders, and so there is no real symmetry 
	between $L_1,L_2$, but the argument for the case when $s\in L_2$ is exactly the same, with $L_1,L_2$ exchanged.) 
	Now 
	$$W_r\cap W_t\subseteq V(G_1)\cap V(G_2) = S,$$
	so it suffices to show that $W_r\cap S\subseteq W_s$. Let $v\in W_r\cap S$. Thus, $v$ is a right-limit vertex in 
	$(L_1, (W_t:t\in L_1))$, and since $r\le_{L_1} s$ it follows that $v\in W_s$. This proves that $W_r\cap W_t\subseteq W_s$,
	and so $(L,(W_t:t\in L))$ is a line-decomposition of $G$. Evidently its width is at most $k'$, and it is a wo-decomposition.
	This proves \ref{concat2}.~\bbox

	We will prove the following, which immediately implies \ref{mainthm} by taking $Z_1,Z_2=\emptyset$:
	\begin{thm}\label{addlimit}
		Let $(L,(W_t:t\in L))$ be a line-decomposition of $G$ from $Z_1$ to $Z_2$, of width at most $k$, and let $m\le k$ such that every split has size 
		at least $m$, and $|Z_1|\le m$. Then $G$ admits a wo-decomposition
		from $Z_1$ to $Z_2$ of width at most $2k-|Z_1|$.
	\end{thm}
	\Proof We proceed by induction on $k$, and for fixed $k$, by induction on $k-m$. As we discussed in the previous section, 
	we may assume that $(L,(W_t:t\in L))$ is tidy. We observe first:
	\\
	\\
	(1) {\em We may assume that $G$ is connected.}
	\\
	\\
	Suppose that $G$ is not connected, and let $C_1$ be a component of $G$. Let $I_1$ be the set of $t\in L$ such that 
	$W_t\cap V(C_1)\ne \emptyset$. Since each $W_t$ is a clique, it follows that $W_t\subseteq V(C_1)$ for each $t\in I_1$. 
	Moreover, we claim that $I_1$ 
	is an interval of $L$. To see this, let $r,t\in I_1$ and suppose that $r\le_L s\le_L t$. There is a path of $C_1$ joining $W_r,W_t$,
	and every such path intersects $W_s$, and so $s\in I_1$. Thus, $I_1$ is an interval.
	Since $I_1\ne L$ (because $G$ is not connected, and so some $W_t$ is disjoint from $V(C_1)$ and hence $t\notin I_1$), 
	we may assume from the symmetry that there exists $s\in L$ such that $s\le_L t$ for all $t\in I_1$, and $s\notin I_1$.
	Let $J$ be the set of $s\in L$ with this property; then $J$ is an initial interval. Let $S$ be the $J$-split. 
	We claim that $S=\emptyset$. Suppose not, and let $v\in S$. Choose $r,s\in L$ such that $v\in W_r\cap W_s$ and 
	$r\in J$ and $s\notin J$. Since $r\in J$, $W_r\cap V(C_1)=\emptyset$; so $W_s\not\subseteq V(C_1)$, and therefore $s\notin I_1$. 
	Since $s\notin J$, and $s\notin I_1$, the definition of $J$ implies that 
	there exists $t\in I_1$ such that $s\not\le_L t$. Thus, $r\le_L t\le_L s$, and $v\in W_r\cap W_s$ and $v\notin W_t$,
	a contradiction. This proves that $S=\emptyset$.  Since $S$ is a split and all splits have size at least $m$, it follows that $m=0$,
	and therefore $Z_1=\emptyset$.

	Since $(L,(W_t:t\in L))$ is tidy, it follows that $Z_2$ is a clique (because it is a subset of a bag of $(L,(W_t:t\in L))$);
	and therefore there is a component $C^*$ of $G$ with $Z_2\subseteq C^*$.  
	We can choose a well-ordering of the set of components of $G$ different from $C^*$; choose an ordinal $\beta$ such that 
	the components of $G$ can be numbered $C_\alpha\;(\alpha<\beta)$, where $\alpha$ ranges over all ordinals $<\beta$. 
	Define $C_{\beta} = C^*$. 
	
	Suppose that the result is true for all components of $G$; thus, for all $\alpha<\beta$, there is a 
	wo-decomposition
	from $\emptyset$ to $\emptyset$ of width at most $2k$; and for $C_{\beta}$ there is a
	wo-decomposition
	from $\emptyset$ to $Z_2$ of width at most $2k$. Then by concatenating these in the natural way, we obtain a wo-decomposition
	of $G$ of width at most $2k-|Z_1|$ (since $Z_1=\emptyset$), as required. This proves (1).
	
	\bigskip
	
	Next, suppose that there are no splits with size $m$. If $m\le k-1$, then the hypotheses all remain true with $m$ replaced by $m+1$,
	and the result follows from the second inductive hypothesis. If $m=k$, then since all splits have size at most $m$ (because
	every split is a subset of a bag), we deduce that there are no splits, and so 
	$|L|=1$; and since $k\le 2k-|Z_1|$, the result is true. Hence we may assume that there is a split of size $m$. Consequently, $m\ge 1$,
	since $G$ is connected and therefore all splits are non-null.
	
	From \ref{countable}, 
	we can number the splits with size exactly $m$ as
	$\{S_i:i\in K\}$, where $K$ is an interval of $\mathbb{Z}$, in such a way that $S_i$ is before $S_j$ for all $i,j\in K$ with $i\le j$. 
	For each $i\in K$, let $\mac A_i$ be the set of all initial intervals $I$ of $L$ such that $S_i$ is the $I$-split.
	Let $I_i= \bigcap_{I\in \mac A_i}I$, and $J_i= \bigcup_{I\in \mac A_i}I$. By \ref{doublesplit}, 
	either $I_i=\emptyset$ and $S_i$ is the set 
	of all left-limit vertices of $(L,(W_t:t\in L))$, or $I_i$ is the minimal element of $\mac A_i$; and either
	$J_i = L$ and $S_i$ is the set of all right-limit vertices, or $J_i$ is the maximal element of $\mac A_i$. 
	\\
	\\
	(2) {\em If $I_i=\emptyset$ then $i$ is the smallest member of $K$,  and if $J_i=L$ then $i$ is the largest member of $K$.}
	\\
	\\
	Suppose that $I_i=\emptyset$. Then $I_i\notin \mac A_i$, since it is not an initial interval of $L$; and so  $S_i$ is the set
	of all left-limit vertices of $(L,(W_t:t\in L))$. Suppose that $i$ is not the smallest member of $K$, and so $i-1\in K$. 
	Choose $H_i\in \mac A_{i}$ and $H_{i-1}\in \mac A_{i-1}$. For each $v\in S_i$, we have $v\in W(L\setminus H_i)\subseteq W(L\setminus H_{i-1})$.
	Moreover, $v\in W_t$ for some $t\in L$, and hence $v\in W_s$ for some $s\in H_{i-1}$ since $v$ is a left-limit vertex.
	Consequently $v\in W(H_{i-1})\cap W(L\setminus H_{i-1}) = S_{i-1}$. This proves that $v\in S_{i-1}$ for each $v\in S_i$,
	which is impossible since $S_{i-1}, S_i$ are distinct and have the same size. Consequently, $i$ is the smallest
	member of $K$. The second statement is proved similarly. This proves (2).
	\\
	\\
	(3) {\em Either $Z_1\subseteq S_i$ for some $i\in K$, or $K$ has a smallest member. Similarly, either $Z_2\subseteq S_i$ for some $i\in K$, or $K$ has a largest member.}
	\\
	\\
	Suppose that there is no $i\in K$ with $Z_1\subseteq S_i$. Since $Z_1$ is finite, and 
	$S_i\cap Z_1\supseteq S_j\cap Z_1$ for all $i,j\in K$ with $i\le j$, it follows that  there exists $z\in Z_1$ such that there is 
	no $i\in K$ with 
	$z\in S_i$. Let $I$ be the intersection of all the initial intervals in
	$\bigcup_{i\in K}\mac A_i$.  Since $z\in Z_1$ is a left-limit vertex, and $z$ belongs to none of the splits $S_i$, it follows that
	$I$ contains all $t\in L$ such that $z\in W_t$; and in particular, $I\ne \emptyset$. Consequently $I$ is an initial interval; let 
	$S$ be the $I$-split. Suppose that there exists $v\in S$ that belongs to none of the splits $S_i\;(i\in K)$. Since $v$ belongs 
	to the $I$-split, there exist $s\in I$ and $t\in L\setminus I$ with $v\in W_s\cap W_t$. Since $t\in L\setminus I$, there exists 
	$i\in K$ and $H\in \mac A_i$ such that $t\in L\setminus H$. But $s\in I\subseteq H$, and so 
	$$v\in W_s\cap W_t\subseteq W(H)\cap W(L\setminus H) = S_i,$$
	a contradiction. Thus, for each $v\in S$ there exists $i\in K$ with $v\in S_i$. Since $S$ is before each $S_i$, it follows that
	$S\cap S_i\supseteq S\cap S_j$ for all $i,j\in K$ with $i\le j$; and so there exists $j\in K$ with $S\subseteq S_j$. But 
	all splits have size at least $m$, and $S$ is a split, and $S_j$ has size $m$; so $S_j=S$. Since $S$ is before $S_i$
	for all $i\in K$, it follows that $j$ is the smallest member of $K$. The second statement of (3) is proved similarly. 
	This proves (3).
	\\
	\\
	(4) {\em We may assume that $K$ has a smallest element.}
	\\
	\\
	Suppose that $K$ has no smallest element. 
	By (3), there exists $i\in K$ such that $Z_1\subseteq S_i$, and so we may assume 
	that $Z_1\subseteq S_0$.
	Let $G_1=G[W(I_0)]$ and $G_2 = G[W(L\setminus I_0)]$. Thus,  $V(G_1)\cap V(G_2) = S_0$. 
	Suppose that:
	\begin{itemize}
		\item Either $Z_1\ne \emptyset$, or $G_1$ admits a wo-decomposition from $S_0$ to $\emptyset$ of width at most $2k-|S_0|$; and
		\item $G_2$ admits a wo-decomposition from $S_0$ to $Z_2$ of width at most $2k-|S_0|$.
	\end{itemize}
	We claim that 
	$G_1$ admits a wo-decomposition from $Z_1$ to $S_0$ of width at most $2k-|Z_1|$.
	To see this, there are two cases. If $Z_1=\emptyset$, then
	adding $S_0$ to every bag of the decomposition given in the first bullet gives
	a wo-decomposition from $S_0$ to $S_0$, and hence from $Z_1$ to $S_0$, as required (since $2k-|Z_1|=2k$).
	Now we
	assume that $Z_1\ne \emptyset$.  By \ref{concat1}, 
	$G_1$ admits a line-decomposition from $Z_1$ to $S_0$ of width at most $k$, and all its bags include $Z_1$, since $Z_1\subseteq S_0$.
	By removing $Z_1$ from all these bags,
	we obtain a  line-decomposition of $G_1\setminus Z_1$ from $\emptyset$ to $S_0\setminus Z_1$ 
	of width at most $k-|Z_1|$. By reversing its order, we obtain a line-decomposition of $G_1\setminus Z_1$ from $S_0\setminus Z_1$ 
	to $\emptyset$
	of width at most $k-|Z_1|$. Since $Z_1\ne \emptyset$, the first inductive hypothesis implies that $G_1\setminus Z_1$
	admits a wo-decomposition from $S_0\setminus Z_1$ to $\emptyset$ of width at most $2(k-|Z_1|) - |S_0\setminus Z_1|$.
	By adding $S_0$
	to all the bags of this wo-decomposition, we have proved that $G_1$ admits a wo-decomposition of width at most
	$$2(k-|Z_1|) - |S_0\setminus Z_1| +|S_0|=2k-|Z_1|$$
	from $S_0$ to $S_0$, and hence from $Z_1$ to $S_0$ since $Z_1\subseteq S_0$, as claimed.
	
	From the second bullet above, $G_2$ admits a wo-decomposition from $S_0$ to $Z_2$ of width at most $2k-|Z_1|$, since $|Z_1|\le |S_0|$.
	From \ref{concat2}, we deduce that $G$ admits a wo-decomposition from $Z_1$ to $Z_2$ of width at most
	$2k-|Z_1|$, so the statement of the theorem holds. 
	
	Consequently, to prove (4), it suffices to prove the two bullets above. But in both cases, we claim that
	the corresponding graph  $G_1$ or $G_2$ admits 
	a line-decomposition $\mac W$ of width at most $k$ (from $S_0$ to $\emptyset$ in the first case, and from $S_0$ to $Z_2$ in the second)
	in which all splits have size at least $m$, and either $\mac W$ has no split of 
	size $m$, or 
	there is one that is before all other splits of $\mac W$ of size $m$. To see this, for $G_1$ we take as $\mac W$ 
	the line-decomposition 
	$(\overleftarrow{I_0},(W_t:t\in \overleftarrow{I_0}))$, where $\overleftarrow{I_0}$ is the line obtained from $I_0$ by reversing its order
	(then all splits of $\mac W$ have size at least $m$, and $S_{-1}$ is a split of $\mac W$ of size $m$, and is before all other splits of size $m$; note that $S_0$ is not  
	a split of $\mac W$). For $G_2$ we take as $\mac W$ the line-decomposition
	$(L\setminus I_0, (W_t:t\in L\setminus I_0))$; then if $1\in K$, $S_1$ is a split of $\mac W$ that is before all other splits of the
	minimum size $m$, and if $1\notin K$ then $\mac W$ has no split of size at most $m$. This proves that $\mac W$ exists with the 
	properties claimed. But if the theorem holds for $\mac W$, then the two bullets above hold, so it suffices to prove the theorem for $\mac W$. This proves (4).

	\bigskip
	
	For each $i\in K$, 
	if $I_i = J_i$, let 
	$G_i$ be the complete graph with vertex set $S_i$. If $I_i\ne J_i$, let 
	$G_i$ be the union of the graphs $G[W_t]$ over all $t\in J_i\setminus I_i$.
	\\
	\\
	(5) {\em For each $i\in K$, there is a wo-decomposition $\mac W_i$ of $G_i$ from $S_i$ to $S_i$ with width at most $2k-m$.}
	\\
	\\
	This is trivial if $I_i=J_i$, so we assume that 
	$J_i\ne I_i$; let $M_i = J_i\setminus I_i$, with order inherited from $L$.
	It follows that $M_i$ is a line,
	and therefore $(M_i,(W_t:t\in M_i))$ is a line-decomposition of $G_i$ from $S_i$ to $S_i$ (because $S_i$ is a subset of $W_t$ for 
	each $t\in M_i$). Thus, $(M_i, (W_t\setminus S_i:t\in M_i))$ is a line-decomposition of $G_i\setminus S_i$ from $\emptyset$ 
	to $\emptyset$, of width at most $k-m$,
	since $S_i$ is a subset of $W_t$ for 
	each $t\in M_i$ and $|S_i|=m$. Since $m\ge 1$, the first inductive hypothesis tells us that there is a wo-decomposition of 
	$G_i\setminus S_i$ from $\emptyset$
	to $\emptyset$, of width at most $2(k-m)$, and by adding $S_i$ to each of its bags, we obtain a wo-decomposition of $G_i$ from 
	$S_i$ to $S_i$, of width at most $2k-m$. This proves (5).
	
	\bigskip
	For each $i\in K$, we define $G'_i$ as follows. If $J_i=L$ (and so $i$ is the largest member of $K$, by (2)),
	let $G_i'$ be the complete graph with vertex set $S_i$.
	If $J_i\ne L$, and $i+1\notin K$, let $G'_i=G[W(L\setminus J_i)]$. If $i+1\in K$, let $G'_i=iG[W(I_{i+1}\setminus J_i)]$.
	\\
	\\
	(6) {\em For each $i\in K$, there is a wo-decomposition $\mac W'_i$
		of $G_i'$ from $S_i$ to $S_{i+1}$ if $i+1\in K$, and from $S_i$ to $Z_2$ if $i+1\notin K$, of width at most $2k-m$.}
	\\
	\\
	First, we assume that $i+1\in K$. Since $S_i$ is before $S_{i+1}$, it follows that $J_i\subseteq I_{i+1}$.
	By (2), $J_i\in \mac A_i$ and $I_{i+1}\in \mac A_{i+1}$, and therefore $I_{i+1}\ne J_i$, since $S_i\ne S_{i+1}$. Consequently
	$I_{i+1}\setminus J_i\ne \emptyset$.
	Let 
	$M_i'= I_{i+1}\setminus J_i$, with order inherited from $L$.
	It follows that $(M_i',(W_t:t\in M_i'))$ is a line-decomposition from $S_i$ to $S_{i+1}$ of $G'_i$
	of width at most $k$; and all its splits have size at least $m+1$. If $m=k$, since all splits of $(M_i,(W_t:t\in M_i'))$ are splits of $(L,(W_t:t\in L))$ and hence have size at most $m$, it follows that $(M_i,(W_t:t\in M_i'))$ has no splits, and hence 
	$|M_i'|=1$. In that case $(M_i',(W_t:t\in M_i'))$ satisfies the statement of (6). Thus, we may assume that 
	$m<k$. From the second inductive hypothesis, applied to $(M_i,(W_t:t\in M_i'))$, since all its splits have size at least $m+1$,
	there is a wo-decomposition
	of $G_i'$ from $S_i$ to $S_{i+1}$ of width at most $2k-|S_i| = 2k-m$. 
	
	Now we assume that $i+1\notin K$, and so $i$ is the largest member of $K$. If $J_i=L$, the claim is clear, using a 
	line-decomposition with only one bag $S_i$; so we assume that $J_i\ne L$. Thus, $G'_i=G[W(L\setminus J_i)]$. Since $J_i\ne L$,
	$(L\setminus J_i,(W_t:t\in L\setminus J_i))$ is a line-decomposition of $G_i$' from $S_i$ to $Z_2$ of width at most $k$, 
	and all its splits have
	size at least $m+1$ (since $i$ is the largest member of $K$). Again, the claim follows from the second inductive hypothesis if $m<k$,
	or trivially if $m=k$. This proves (6).
	
	\bigskip
	
	From (4), we may assume that $1$ is the smallest member of $K$.
	Define $G_0' = G[W(I_1)]$ if $I_1\ne \emptyset$, and let $G_0'$ be the complete graph with vertex set $S_1$ if $I_1=\emptyset$. 
	Then, similarly, there is a wo-decomposition $\mac W'_0$ of $G_0'$ from $Z_1$ to $S_1$ of width at most $2k-|Z_1|$.
	By \ref{concat2} and concatenating
	$$\mac W'_0, \mac W_1, \mac W_1'\LL \mac W_{i-1}', \mac W_i, \mac W_i'$$
	if $K$ has a largest member $i$, or by concatenating 
	$$\mac W'_0, \mac W_1, \mac W_1', \mac W_2,\ldots$$
	otherwise, we deduce that $G$ admits a wo-decomposition of width at most $2k-|Z_1|$ from $Z_1$ to $Z_2$. 
	(Note that in the case when $K$ has no largest member, $Z_2$ is included in $S_i$ for all sufficiently large $i$, by (3), and so
	$Z_2$ is a set of right-limit vertices of the wo-decomposition we just constructed.) This proves \ref{addlimit}.~\bbox
	
	\section{Prime decompositions}
	There is another route we could take to making line-decompositions as palatable as possible. In the proof of step (1) of
	\ref{addlimit}, we used ordinals of arbitrary magnitude, but all the remainder of the proof of \ref{addlimit} was just concatenating sequences
	of wo-decompositions to make longer wo-decompositions, and this process was only iterated $k$ times. If we could find a way to
	avoid step (1), we might get something nicer.
	
	That motivates looking at ``prime'' line-decompositions.
	Let us say a
	line-decomposition
	$(L,(W_t:t\in L))$ is {\em prime} if it is tidy, and $G$ is connected, and the $I$-split is different from the $J$-split, for all distinct initial intervals $I,J$.  If we run the proof of \ref{addlimit} starting with a prime line-decomposition, three nice things happen:
	\begin{itemize}
		\item $G$ is already connected, so we do not need step (1) of the proof;
		\item the graphs $G_i$ used in step (5) of the proof are all trivial, because $I_i=J_i$ in the notation of that proof, so we could skip step (5);
		\item the line-decompositions $(M_i,(W_t:t\in M_i'))$ of the graphs $G_i'$ used in step (6) are also prime, so we can apply the inductive hypothesis to them.
	\end{itemize}
	We deduce:
	\begin{thm}\label{usingprime}
		If $(L,(W_t:t\in L))$ is a prime line-decomposition of $G$, of width at most $k$, then $L$ has order-type that of a subset of $\mathbb{Z}^k$,
		and $G$ admits a wo-decomposition $(L',(W_t':t\in L'))$ of width at most $2k$ where $L'$ is an ordinal and $L'\le \omega^k$.
	\end{thm}
	
	Of course, not every connected graph that admits a tidy line-decomposition admits a prime line-decomposition of the same width. For instance,
	a {\em caterpillar} is a tree in which some path passes through every vertex of degree more than two. A connected graph $G$ admits
	a
	tidy line-decomposition of width one if and only if $G$ is a caterpillar, but $G$ admits a prime line-decomposition of width one
	if and only if $G$ is a path. Indeed, from the first statement of \ref{usingprime}, if $G$ admits a prime line-decomposition
	of width at most $k$, then $V(G)$ is countable.
	
	On the other hand, we will show that a general tidy line-decomposition can be built
	by substituting prime line-decompositions in one another (iterating only $k$ times, where $k$ is the width).
	Let us say this
	more exactly. Let $L$ be a linear order, and let $\mac I$ be a set of some of its initial intervals. For each $I\in \mac I$,
	let $L_I$ be a linear order, such that $L$ and all the orders
	$L_I$ are pairwise disjoint. Let $M$ be the linear order with element set $L\cup \bigcup_{I\in \mac I}L_I$, where $s\le t$ if either:
	\begin{itemize}
		\item for some $I\in \mac I$, $s,t\in L_I$, and $s\le t$ in the order of $L_I$; or
		\item $s,t\in L$, and $s\le t$ in the order of $L$; or
		\item $s\in L$, and $t\in L_I$ for some $I\in \mac I$, and $s\in I$; or
		\item $s\in L_I$ for $I\in \mac I$, and $t\in L$, and $t\notin I$.
	\end{itemize}
	It is easy to check that $M$ is indeed a linear order. Now let $(L, (W_t:t\in L))$ be a line-decomposition of some graph $G$, and for each
	$I\in \mac I$ let $(L_I, (W^I_t:t\in L_I))$ be a line-decomposition of some non-null graph $G_I$, where the graphs $G$ and $G_I\;(I\in \mac I)$
	are all pairwise vertex-disjoint. Let $H$ be the graph obtained from the disjoint union of $G$ and the graphs $G_I\;(I\in \mac I)$
	by adding an edge between $u,v$ for each $I\in \mac I$, each $u\in V(G_I)$ and each $v\in W(I)\cap W(L\setminus I)$.
	Define $M$ as before, and for each $t\in M$, define $X_t=W_t$ if $t\in L$, and $X_t=W_t\cup (W(I)\cap W(L\setminus I))$ if $t\in L_I$
	for some $I\in \mac I$. Then $(M,(X_t:t\in M))$ is a line-decomposition of $H$, tidy if $(L, (W_t:t\in L))$ and each
	$(L_I, (W^I_t:t\in L_I))$ are tidy; and we say that $(M,(X_t:t\in M))$ is obtained by
	{\em substituting $((L_I, (W^I_t:t\in L_I))\;I\in \mac I)$ into
		$(L, (W_t:t\in L))$}. We will sketch a proof that:
	\begin{thm}\label{subst}
		Let $(M,(X_t:t\in M))$ be a tidy line-decomposition of some connected graph $H$. Then there exist $G,(L, (W_t:t\in L))$, $\mac I$,
		$G_I\;(I\in \mac I)$, and $(L_I, (W^I_t:t\in L_I))$ for each $I\in \mac I$, as above, such that
		\begin{itemize}
			\item $(M,(X_t:t\in M))$ is obtained by substituting $((L_I, (W^I_t:t\in L_I))\;I\in \mac I)$ into
			$(L, (W_t:t\in L))$;
			\item $(L, (W_t:t\in L))$ is prime, with width at most that of $(M,(X_t:t\in M))$; and
			\item for each $I\in \mac I$, $(L_I, (W^I_t:t\in L_I))$ is tidy and has width strictly smaller than that of $(M,(X_t:t\in M))$.
		\end{itemize}
	\end{thm}
	\Proof For each split $C$ of $(M,(X_t:t\in M))$, let $F_C$ be the set of all initial intervals $I$ of $M$ such that the
	$I$-split of $(M,(X_t:t\in M))$ equals $C$. If $|F_C|>1$, we say $C$ is a {\em repeated split}. If $C$ is a repeated split,
	the set of all $t\in M$ such that $t$ belongs to some but not all members of $F_C$ is an interval of $M$; let us call this interval
	the {\em repeat interval of $C$} and denote it by $L_C$.
	\\
	\\
	(1) {\em If $C$ is a repeated split of $(M,(X_t:t\in M))$ and $t\in L_C$, then $C\subseteq X_t$.}
	\\
	\\
	Since $t\in L_C$, there exist initial intervals $I,J\in F_C$ such that $t\in J\setminus I$. Let $v\in C$.
	Thus, $v$ belongs to the $I$-split, and therefore to $X(I)$; and so there exists $s\in I$ with $v\in X_s$. Similarly there exists $s'\in M\setminus J$ such that $v\in X_{s'}$. Since $s\le t\le s'$, it follows that $v\in X_t$. This proves (1).
	\\
	\\
	(2) {\em If $C$ is a repeated split of $(M,(X_t:t\in M))$, and $C'$ is a split with $C\not\subseteq C'$,
		and $J\in F_{C'}$, then either $L_C\cap J=\emptyset$ or $L_C\subseteq J$.}
	\\
	\\
	Choose $v\in C\setminus C'$. Since $v\notin C'=X(J)\cap X(M\setminus J)$, either $v\notin X(J)$ or $v\notin X(M\setminus J)$.
	But $v\in X_t$ for all $t\in I$, by (1); and so either $L_C\cap J=\emptyset$ or $L_C\cap (M\setminus J)=\emptyset$. This proves (2).
	\\
	\\
	(3) {\em If $C,C'$ are distinct repeated splits of $(M,(X_t:t\in M))$, then either $L_C\subseteq L_{C'}$
		and $C'$ is a proper subset of $C$,
		or $L_{C'}\subseteq L_C$ and $C$ is a proper subset of $C'$, or $L_C\cap L_{C'}=\emptyset$.}
	\\
	\\
	We may assume that $C\ne C'$, and so without loss of generality, we assume that $C\not\subseteq C'$. By (2), for each initial interval
	$J\in F_{C'}$, either  $L_C\cap J=\emptyset$ or $L_C\subseteq J$. If $L_C\cap J=\emptyset$ for every $J\in F_{C'}$, then
	$L_C\cap L_{C'}=\emptyset$; and
	if $L_C\subseteq J$ for all $J\in F_{C'}$ then again $L_C\cap L_{C'}=\emptyset$. Finally, we assume that there exist $J_1,J_2\in F_{C'}$
	such that $L_C\cap J_1=\emptyset$ and  $L_C\subseteq J_2$. It follows that $L_C\subseteq L_{C'}$. If $C'\subseteq C$,
	then $C'$ is a proper subset of $C$ and the claim holds, so we may assume that $C'\not\subseteq C$.
	By the same argument with $C,C'$ exchanged,
	we deduce that either
	\begin{itemize}
		\item $L_C\cap L_{C'}=\emptyset$ (and this is false since $\emptyset\ne L_C\subseteq L_{C'}$); or
		\item $L_{C'}\subseteq L_C$, and so $L_{C'}= L_C$.
	\end{itemize}
	Thus, it remains to handle the case when $L_{C'}= L_C$ and neither of $C,C'$ is a subset of the other. Let $J\in F_C$, and suppose that
	$L_C \cap J\ne \emptyset$ and $L_C\not\subseteq J$. Then $J$ includes some but not all of $L_C$, and so there exist $s,t\in L_C$
	with $s\in J$ and $t\notin J$. Hence $X_s\cap X_t$ is a subset of the $J$-split, that is, of $C$; but by (1), $X_s\cap X_t$
	includes $C'$ since $s,t\in L_{C'}=L_C$, so
	$C'\subseteq C$, a contradiction. Thus, if $J\in F_C$ then either $L_C \cap J= \emptyset$ or $L_C\subseteq J$.
	Let $I_1$ be the set of all
	$s\in M$ such that $s\notin L_C$ and $s<t$ for all $t\in L_C$. Thus, $I_1$ is either an initial interval or empty. Similarly,
	let $I_2$ be the set of all $s\in M$ such that $s\le t$ for some $t\in L_C$; then $I_2$ is either an initial interval or equals $M$.
	
	We claim that if $J\in F_C$ then $J=I_1$ or $J=I_2$. To see this, choose $s\in L_C$. Suppose first that $L_C\cap J=\emptyset$,
	so $J\subseteq I_1$.
	If there exists $t\in I_1\setminus J$, then $t$ belongs to some but not all of the members of $F_C$, because it does not belong to $J$,
	and belongs to any member of $F_C$ that contains $s$ (and there is such a member, by definition of $L_C$). But then $t\in L_C$, a
	contradiction. So if $L_C\cap J=\emptyset$ then $J=I_1$. Next suppose that $L_C\subseteq J$, If there exists $t\in J\setminus I_2$,
	then $t$ belongs to some but not all of the members of $F_C$, because it belongs to $J$,
	and does not belong to any member of $F_C$ that does not contains $s$ (and there is such a member, by definition of $L_C$).
	But then again $t\in L_C$, a contradiction. This proves the claim that if $J\in F_C$ then either $J=I_1$ or $J=I_2$. Since $|F_C|\ge 2$,
	it follows that $F_C = \{I_1,I_2\}$. In particular, $I_1,I_2$ are both initial intervals. Similarly $F_{C'} = \{I_1,I_2\}$; but then the $I_1$-split equals both $C$ and $C'$,
	a contradiction. This proves (3).
	
	\bigskip
	
	Let $\mac C$ be the set of all repeated splits $C$ such that $L_C$ is maximal. By (3), the sets $L_C\;(C\in \mac C)$
	are disjoint. Moreover:
	\\
	\\
	(4) {\em For every repeated split $C$, there exists $C'\in \mac C$ such that $L_C\subseteq L_{C'}$.}
	\\
	\\
	Choose a repeated split $C'$ with $C'$ minimal such that $L_C\subseteq L_{C'}$.
	Then by (3), $C'\in \mac C$. This proves (4).
	
	\bigskip
	
	Let $L$ be the set of all $t\in M$ that do not belong to $\bigcup_{C\in \mac C} L_C$, with the linear order inherited from that of $M$.
	For each $t\in L$, let $W_t=X_t$; then
	$(L,(W_t:t\in L))$ is a line-decomposition of some subgraph $G$ of $H$. For each $C\in \mac C$, and each $t\in L_C$, define
	$W^C_t=X_t\setminus C$; then $(L_C, (W^C_t:t\in L_C))$ is a tidy line-decomposition of some subgraph $G_C$ of $H$, of width at most the width of
	$(M,(X_t:t\in M))$ minus $|C|$. Since each $C\ne \emptyset$ (because $H$ is connected), it follows that each
	$(L_C, (W^C_t:t\in L_C))$ has width less than that of $(M,(X_t:t\in M))$. To complete the proof, we need to check that 
	$(L,(W_t:t\in L))$ is prime. If $I$ is an initial interval of $L$, then there is an initial interval $J$ of $M$ such that
	the split of $I$ in $(L,(W_t:t\in L))$ equals the split of $J$ in $(M,(X_t:t\in M))$; and if $I\subseteq I'$ are distinct initial intervals of
	$L$, then
	the corresponding initial intervals $J\subseteq J'$ of $M$ satisfy that $J'\setminus J$ is not a subset of any member of $\mac C$,
	since $I' \setminus I \subseteq J' \setminus J$.
	By (4), the splits of $J$ and of $J'$ in $(M,(X_t:t\in M))$ are distinct, and so the splits of $I$ and of $I'$ in $(L,(W_t:t\in L))$
	are distinct. Since $H$ is connected, and is obtained from the disjoint union of the graphs $G$ and the graphs $G_C\;(C\in \mac C)$ by adding edges from each $G_C$ to the clique $C$ of $G$, it follows that $G$ is connected. Thus, $(L,(W_t:t\in L))$ is prime. This proves \ref{subst}.~\bbox
	
	This result shows that we can build a tidy line-decomposition of width $k$ of a {\em connected} graph by starting with a prime
	decomposition and substituting tidy line-decompositions of width $<k$ into it. But the tidy
	line-decomposition we are substituting are decompositions of subgraphs that might not be connected. So to complete a recursive
	structure theorem, we need another
	step, a proof that we can build any tidy line-decomposition of width $\le k$ of a graph $G$ from tidy line-decompositions
	of width $\le k$ of
	the components of $G$, by concatenating the latter in some order. This is not hard and we omit it.

	\section*{Acknowledgements}
	We would like to thank the anonymous referees for helpful comments and suggestions.
	
	\section*{Declaration of generative AI use}
	We acknowledge the use of ChatGPT-5.6 Sol in proofreading this manuscript.
	
	\section*{Remarks on open access}
	For the purpose of open access, the author has applied a CC BY public copyright licence to
	any author accepted manuscript arising from this submission.
	

	\section*{Statements and Declarations}
	\subsection*{Funding}
	Tung Nguyen was supported by a Titchmarsh Research Fellowship and a Christ Church Research Centre Grant.
	Alex Scott was supported by EPSRC grant EP/X013642/1.
	Paul Seymour was supported by AFOSR grant FA9550-22-1-0234, and by NSF grant DMS-21454169.
	\subsection*{Competing interests}
	The authors have no relevant financial or non-financial interests to disclose.
\end{document}